\documentclass[12pt]{amsart}
\usepackage{amsfonts, amsmath, amssymb, amscd, latexsym, graphicx}

\newtheorem{thm}{Theorem}[section]
\newtheorem{prop}[thm]{Proposition}

\newtheorem{defn}[thm]{Definition}

\def\zz{\mathbb{Z}}
\def\rr{\mathbb{R}}

\def\ii{^{-1}}

\renewcommand\>{\rangle}

\def\plr{PL_2(\rr)}
\def\autf{{\rm Aut}\,F}
\def\plpr{PL_2^+(\rr)}
\def\autpf{{\rm Aut}^+F}

\begin{document}

\title[The automorphism group of Thompson's group $F$]{The automorphism group of Thompson's group $F$: subgroups and metric properties}

\author{Jos\'e Burillo}
\address{
Departament de Matem\`atica Aplicada IV, Escola Polit{\`e}cnica
Superior de Castelldefels, Universitat Polit\`ecnica de Catalunya, C/Esteve Torrades 5,
08860 Castelldefels, Barcelona, Spain} \email{burillo@ma4.upc.edu}
\thanks{
The first author acknowledges support from MEC grant MTM2008-01550.}

\author{Sean Cleary}
\address{Department of Mathematics,
The City College of New York \& The CUNY Graduate Center, New York,
NY 10031} \email{cleary@sci.ccny.cuny.edu}
\thanks{The second author  acknowledges support  by the National Science Foundation under Grant No. \#0811002.}

\thanks{}

\begin{abstract}
We describe some of the geometric properties of the automorphism group $Aut(F)$ of Thompson's group $F$.  We give realizations of $Aut(F)$ geometrically via periodic tree pair diagrams, which lead to natural presentations and give effective methods for estimating the word length of elements.  We study some natural subgroups of $Aut(F)$ and their metric properties.  In particular, we show that the subgroup of inner automorphisms of $F$ is at least quadratically distorted in $Aut(F)$, whereas other subgroups of $Aut(F)$ isomorphic to $F$ are undistorted.

\end{abstract}

\maketitle

\section*{Introduction}

Thompson's group $F$ has a number of perplexing properties.  These include a wide range of unusual geometric properties, both on the large scale and the local scale.  These have been studied by a number of authors.
Here, we develop a geometric understanding of the automorphism group of $F$.  The algebraic structure of $Aut(F)$ was described by Brin \cite{brin} via a short exact sequence involving the subgroup of inner automorphisms of $F$ together with the product of two copies of Thompson's group $T$.  The geometric structure of $Aut(F)$ is related to the geometric description of the commensurator of $F$, developed in Burillo, Cleary and R\"over \cite{bcr}.  Here, we describe $Aut(F)$ in terms of periodic infinite tree pair diagrams.  These descriptions lead to a natural presentation for $Aut(F)$ as well as some effective estimates for word length in $Aut(F)$ in terms of the complexity of their eventually periodic tree pair diagrams.  We estimate the word length as a function of the size of the periodic part of the diagrams, together with the size of a remaining non-periodic part. 

There are several natural subgroups to consider of $Aut(F)$ arising naturally via Brin's algebraic description  as a short exact sequence.  We consider several such subgroups, giving geometric descriptions of them and their presentations.  We estimate the word length in these subgroups and compute the exact distortion to be quadratic in $Aut(F)$ of the simplest one of these subgroups.  The construction shows that $F$ itself is quadratically distorted in the most natural of ways in which it lies in $Aut(F)$, as the subgroup of inner automorphisms.  In contrast, there are many other ways in which $F$ lies as a subgroup in $Aut(F)$ and in some of those cases, we show that it is undistorted.

\section{Background on $F$ and its automorphism group}

Thompson's group $F$ is typically defined as a group of homeomorphisms of the unit interval. Elements of $F$ are maps
$$
f:[0,1]\longrightarrow [0,1]
$$
which satisfy
\begin{enumerate}
\item they are orientation-preserving homeomorphisms,
\item they are piecewise-linear, with finitely many breakpoints (points where the linear slope changes),
\item all breakpoints are dyadic integers (lying in $[0,1]\cap \zz[\frac12]$),
\item all slopes of the linear parts are powers of 2.
\end{enumerate}
For an introduction to $F$, and proofs of its most important properties, see Cannon, Floyd and Parry \cite{cfp}.

Thompson's group $T$ is defined in a completely analogous way as $F$, but the homeomorphisms are defined on the circle $S^1$ instead of the unit interval. We usually interpret elements of $T$ also as maps of the interval but with the endpoints identified, and then we can also think of them as piecewise-linear, having dyadic breakpoints, etc. Hence, $F$ is a subgroup of $T$. Again, see \cite{cfp} for details.

\begin{figure}
\centerline{\includegraphics[width=5in]{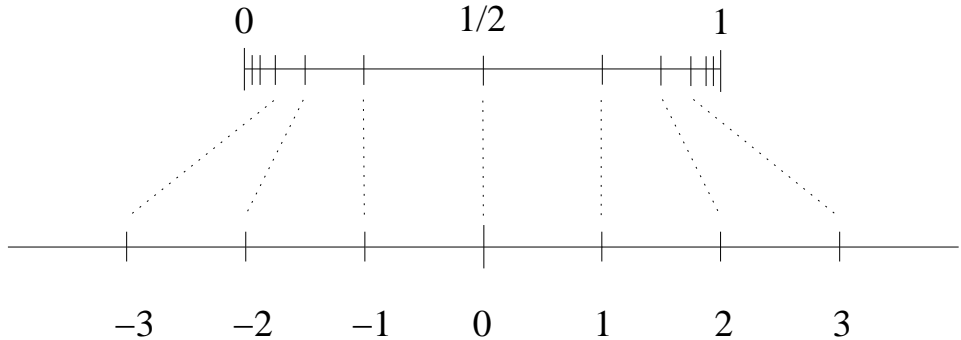}}
\caption{The map used to transport $F$ from $[0,1]$ to the real line.}
\label{conjmap}
\end{figure}

To study the group of automorphisms of $F$, it is convenient to introduce the group $\plr$. This group is the group of piecewise-linear homeomorphisms of the real line which have also dyadic breakpoints and power-of-two slopes, but whose set of breakpoints can be infinite, but discrete. Thompson's group $F$ can be seen as a group of homeomorphisms of the real line as well, instead of $[0,1]$, by conjugating elements with the suitable map (see Figure \ref{conjmap}). Even though it seems that this conjugation could introduce an infinite number of breakpoints, it is easy to see that it does not. The slope $2^k$ near 0, for instance, results in translation by $k$ near $-\infty$ when conjugated into the real line. Hence, $F$ is seen as a subgroup of $\plr$, as those elements which satisfy the following property: there exists a real number $M>0$, and two integers $k,l$ such that
\begin{enumerate}
\item for all $x>M$, we have $f(x)=x+k$
\item for all $x<-M$, we have $f(x)=x+l$
\end{enumerate}
In particular, this implies that all breakpoints  are within the interval $[-M,M]$ and hence the set of breakpoints is finite.

The group $\plr$ is the group where all the action will be contained. Brin \cite{brin}, shows that the group $\plr$ acts as a natural setting for these studies, showing that the automorphism groups of  both Thompson's group $F$ and its commutator subgroup are subgroups of $\plr$:\begin{thm} (Brin \cite{brin}) For a subgroup $G\subset\plr$, let $N(G)$ be its normalizer in $\plr$, that is, all those elements in $\plr$ which conjugate $G$ to itself. Then we have:
\begin{enumerate}
\item $\autf\cong N(F),$
\item $\autf'\cong N(F')=\plr.$
\end{enumerate}
\end{thm}
This theorem is what allows us to investigate the automorphism group of $F$ as  a subgroup of $\plr$, where a piecewise linear map acts on $F$ by conjugation. That is, the group of automorphisms of $F$ is isomorphic to those elements in $\plr$ which conjugate $F$ to itself.

We note that $\plr$ has an index two subgroup, which we will denote $\plpr$, of those maps in $\plr$ which preserve the orientation.  A natural representative of the reversing coset is the map $x\mapsto -x$. The same way, the group of automorphisms of $F$ has also a subgroup $\autpf$, also preserving orientation. Since in both cases the subgroup has index two, we will restrict our study to the orientation-preserving automorphisms, noting that $F$ is actually included in $\plpr$.

The key characteristic of the elements of $F$ inside $\plpr$ is the fact that near infinity, they are translations by integers. The two integers $k$ and $l$ which exist for each element of $F$ give precisely the abelianization map for $F$ (see  Cannon, Floyd and Parry \cite{cfp}). Then, the commutator $F'$ is given by the elements which are the identity near infinity. This shows that $F'$ is exactly the subgroup of $\plpr$ of those elements which have bounded support; that is, those which are the identity outside a bounded interval. This gives an idea of why the automorphism group of $F'$ is the whole $\plr$. See Brin's paper \cite{brin} for details.

To see which elements of $\plpr$ conjugate $F$ to itself, we note that the condition where $f(x)=x+k$, for $x>M$, and the fact that all integers appear somewhere, forces a conjugating element $\alpha$ to satisfy $\alpha(x+1)=\alpha(x)+1$, also for sufficiently large $x$. The counterpart in $-\infty$ also means that this same condition is satisfied on the other side. Hence an element $\alpha\in\autpf$ must satisfy $\alpha(x+1)=\alpha(x)+1$ outside of some bounded interval.

Finally, if we consider a map $\alpha$ on $\rr$ which satisfies $\alpha(x+1)=\alpha(x)+1$, it descends into $S^1$ seen as $\rr/\zz$, so it corresponds to a map of $S^1$. This means that to an element of $\autpf$ we can assign two elements of $T$, one given by the behavior near $\infty$ and another one near $-\infty$ (which have no interaction among them, and hence, commute). This fact is summarized in the following structure theorem for $\autpf$, proved by Brin in \cite{brin}.

\begin{thm}(Brin  \cite{brin})
\begin{enumerate}
\item
The group $\autpf$ is exactly the subgroup of $\plpr$ given by those elements $\alpha$ for which there exists $M>0$ such that $\alpha(x+1)=\alpha(x)+1$, for all $x\notin[-M,M]$.
\item We have a short exact sequence
$$
1\longrightarrow F\longrightarrow \autpf\longrightarrow T\times T\longrightarrow 1.
$$
\end{enumerate}
\end{thm}

In this paper, we introduce an interpretation for $\autpf$ given by binary trees analogous to the one for $F$,  use it to construct a presentation for it, describe some interesting subgroups for $\autpf$, and study the large-scale metric properties of some of these subgroups.

\section{Binary trees}

\begin{figure}
\centerline{\includegraphics{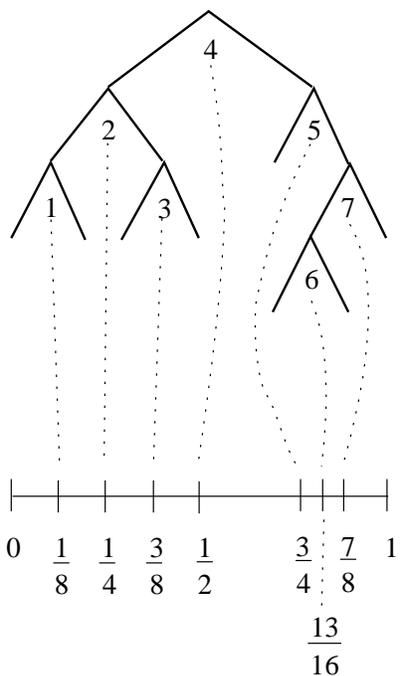}}
\caption{An example of a binary tree and the subdivision it represents in the unit interval.}
\label{caretorder}
\end{figure}

Elements of Thompson's group $F$ can be represented using a pair of binary trees to encode the subdivisions of the unit interval needed to construct the piecewise-linear map in a standard way.  By ``caret'' of a rooted binary tree, we mean an internal node together with its two downward directed edges.  Each caret represents a subdivision of an appropriate subinterval into two halves, which are its children. That way, a binary tree represents a subdivision of the interval. Two trees with the same number of leaves represent an element of $F$, by mapping the leaf subintervals linearly in an order-preserving way. The carets are naturally ordered from left to right, according to the order of the subintervals they represent in $[0,1]$. See Figure \ref{caretorder} for an example.

\begin{figure}
\centerline{\includegraphics[width=5in]{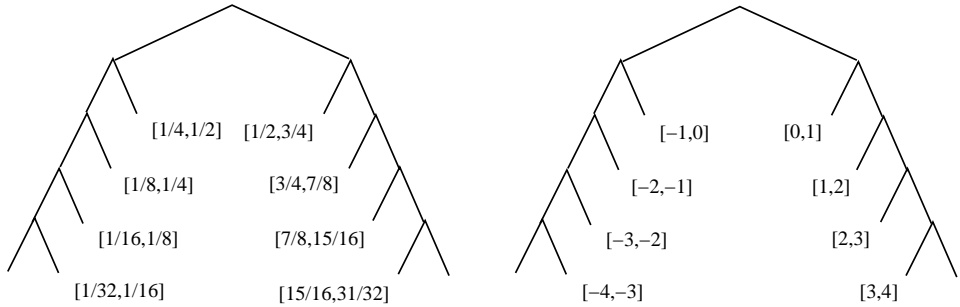}}
\caption{The conjugation between the unit interval and the real line, seen with trees representing elementary subdivisions.}
\label{realtree}
\end{figure}

Our goal is to extend this representation to elements of $\autf$, and even to the whole $\plr$. Since according to the previous section, elements of $\autf$ are represented by piecewise-linear maps on the line instead of $[0,1]$, we need to clarify how pairs of binary trees will correspond to maps of the real line. If we consider the conjugation map which relates the interval to the line, we see how  to construct this representation. The carets located on the right and left branches of a tree, represent the main intervals which appear in the conjugating map, see Figure \ref{conjmap}. Hence, these carets will represent integer intervals of length one in the real line representation, with the exception of the leftmost and rightmost carets which represent intervals to $-\infty$ and $+\infty$ respectively. See Figure \ref{realtree}.

\begin{figure}
\centerline{\includegraphics[width=5in]{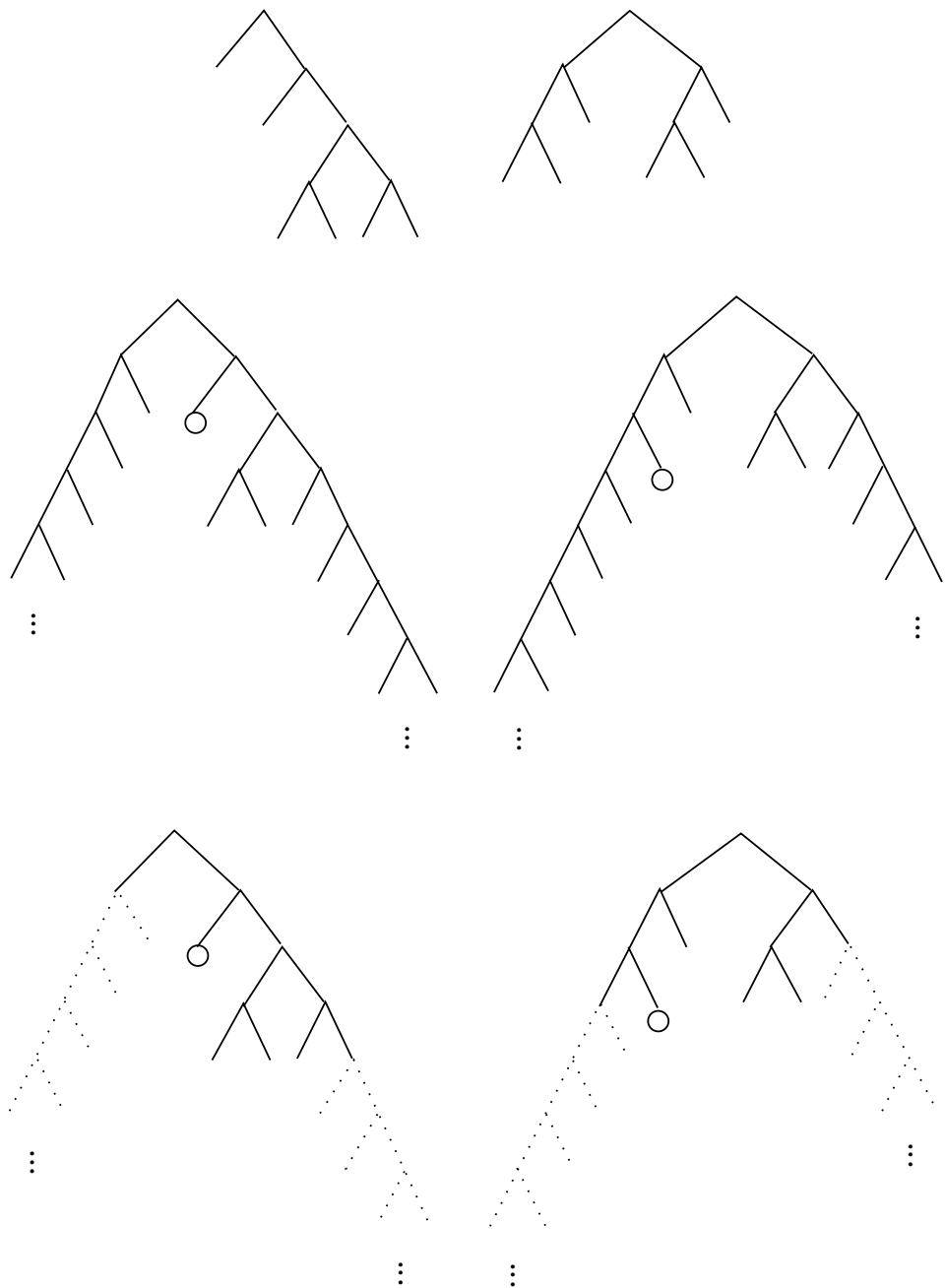}}
\caption{An element of $F$ but now seen as a pair of trees encoding the map of the real line. Here we see: the element in $[0,1]$, then the same element in $\rr$, and in the third picture the dotted carets represent the tails that have been added.}
\label{realelement}
\end{figure}

So we are now considering infinite binary trees (with infinitely many leaves). We need a mark on the two trees to signal a reference starting point. In this way, each element of $\plr$ can be represented by a pair of (possibly infinite) binary trees with a marking. Elements of $F$ correspond to those for which only finitely many of the unit integer intervals are subdivided, and hence the trees are finite with the
exception of  the right and left ``tails." See Figure \ref{realelement} for an example.  The choice of marking is arbitrary; there are many marked diagrams representing the same element.  We note that it is possible to represent the whole $\autf$ by indicating some orientations with the markings.  However, here we will restrict our consideration to the orientation-preserving subgroup $\autpf$ of index two for simplicity.

In the particular case for $\autf$, the pairs of trees have to be periodic. The condition $\alpha(x+1)=\alpha(x)+1$ means that whatever subdivision we have in the interval $[k,k+1]$, for $k\in\zz$ sufficiently large, it must repeat in the following intervals. Hence, the infinite trees representing elements of $\autf$ can be coded in a finite way, because from some point on, the trees are periodic. This  leads to the following definition.

\begin{defn} An \emph{eventually periodic binary tree} is a subtree of the infinite binary tree, for which there exists an integer $N>0$ such that all intervals $[k,k+1]$, for $k\ge N$ have the same subdivision; and also, all intervals $[-k,-k+1]$, for $k\ge N$ also have the same subdivision. The subdivision for the intervals near infinity (the same for all) is not necessarily the same as the one for the intervals near $-\infty$.
\end{defn}

\begin{figure}
  \includegraphics[width=5in]{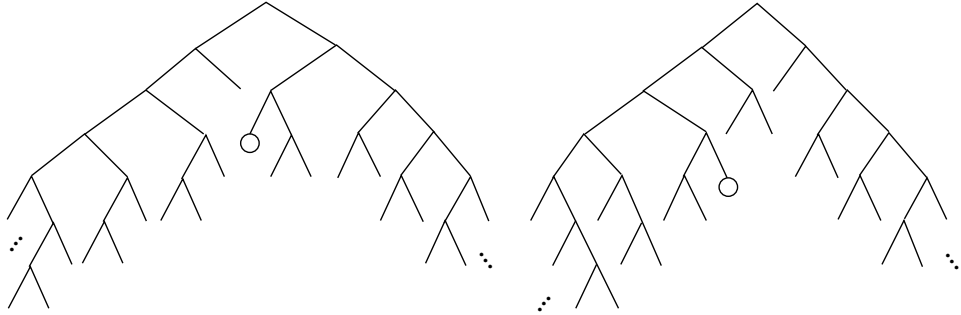}\\
  \caption{An example of an eventually periodic map represented by infinite, also eventually periodic trees.}\label{evperiodic}
\end{figure}

See Figure \ref{evperiodic} for an example of an eventually periodic element. The trees are infinite, but only at the two extremes, and that each integer interval is subdivided only finitely many times. Hence the breakpoints form a discrete subset. And the fact that the trees are eventually periodic means that for each tree, the right and left tails are periodic. Since now the trees are infinite, we need a marking to indicate how the leaves are mapped to each other. The marking (unoriented in the case of $\autpf$) is given by a little circle in a leaf in each tree.  Then, those two leaves are mapped to each other, and the corresponding leaves are mapped in order-preserving way.

It is straightforward to verify that conjugating an element of $F$ by an eventually periodic tree, the result is another element of $F$. The tails get cancelled because they are identical and give redundant subdivisions, capable of being reduced.
%The reader can try conjugating the previous element of $F$ by this eventually periodic element, to be convinced that these conjugations really yield elements of $F$.

\begin{figure}
  \includegraphics[width=12cm]{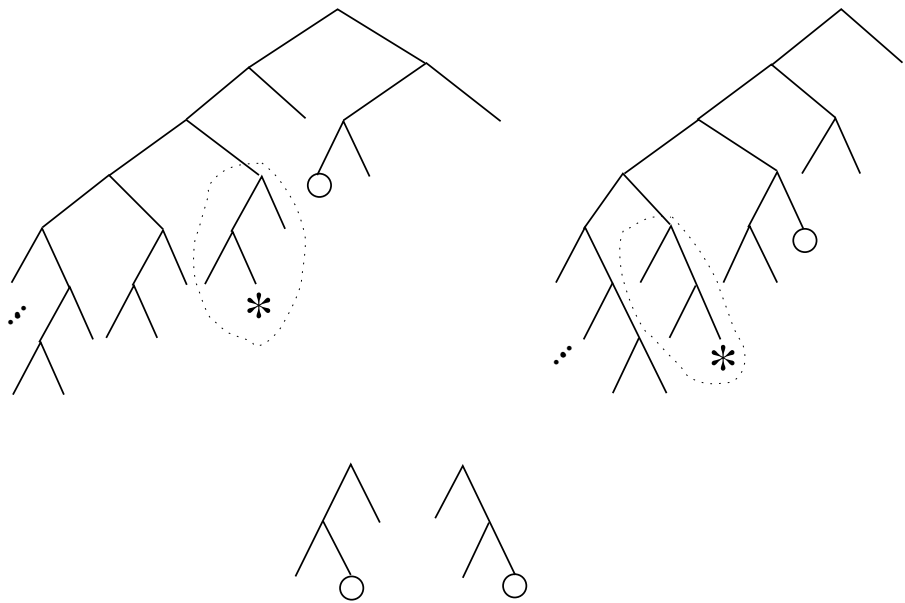}\\
  \caption{The left tail of the element in Figure \ref{evperiodic}. The two trees enclosed with a dotted line are the ones repeating indefinitely down the tail, and the two leaves marked with an asterisk map to each other (they are both at distance three from the marking). The resulting element of $T$ is constructed with these two trees, and is represented on the bottom.}\label{tcrosst}
\end{figure}

Each orientation-preserving automorphism of $F$ can be represented by an eventually periodic map of $\rr$, and hence by a pair of eventually periodic trees. The map from $\autpf$ to $T\times T$ can be read off directly from the infinite tree. The map is given by the trees which subdivide in both tails. See Figure \ref{tcrosst} for the element of $T$ corresponding to the left tail for the element of Figure \ref{evperiodic}.

We note that we are restriction to the binary case; the generalizations $F(p)$ of Thompson's group $F$ have much more complicated automorphism groups; see Brin and Guzman \cite{bringuzman}.

\section{Presentation}

The short exact sequence
$$
1\longrightarrow F\longrightarrow \autpf\longrightarrow T\times T\longrightarrow 1.
$$
can be used to compute a presentation for $\autpf$. We can use the standard presentation for $F$ given by
$$
\<x_0,x_1\,|\,[x_0x_1\ii ,x_0\ii x_1x_0],[x_0x_1\ii ,x_0^{-2}x_1x_0^2]\>.
$$
A standard finite presentation for $T$,  given in Cannon, Floyd and Parry \cite{cfp}, supplements the two generators $x_0$ and $x_1$ for $F$ with a new generator $c$, which is an element of order 3. In our case, we will choose a different torsion element, denoted by $t$, which is diagrammatically smaller and easier to lift to $\autpf$ and is also used in Burillo, Cleary, Stein and Taback \cite{thompt}. The two elements $c$ and $t$ are depicted in Figure \ref{tgens}, and note that they satisfy the relation $c=tx_0$. Using this relation and the presentation for $T$ given in \cite{cfp}, it is easy to see that $T$ is generated by $x_0,x_1$ and $t$ and with the following relators, which are labelled for future reference:
\begin{figure}
  \includegraphics{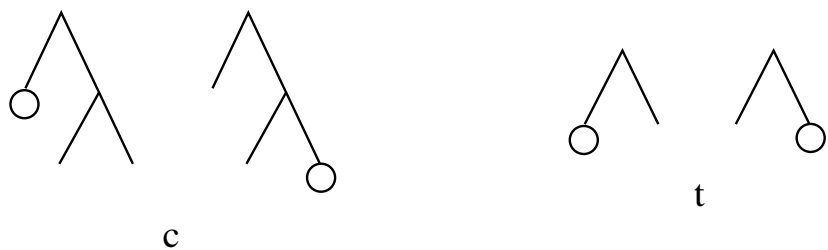}\\
  \caption{Two torsion elements of $T$. The element $c$ is used as a generator in \cite{cfp}, and the element $t$ is used here.}\label{tgens}
\end{figure}
\begin{enumerate}
\item[(f1)] $[x_0x_1\ii ,x_0\ii x_1x_0]$
\item[(f2)] $[x_0x_1\ii ,x_0^{-2}x_1x_0^2]$
\item[(t1)] $x_1\ii tx_0^2x_1\ii x_0\ii t$
\item[(t2)] $x_1\ii tx_0^2x_1\ii tx_0tx_0$
\item[(t3)] $x_1^{-2}tx_0^2x_1\ii x_0\ii t$
\item[(t4)] $t^2$.
\end{enumerate}

\begin{figure}
  \includegraphics[width=5in]{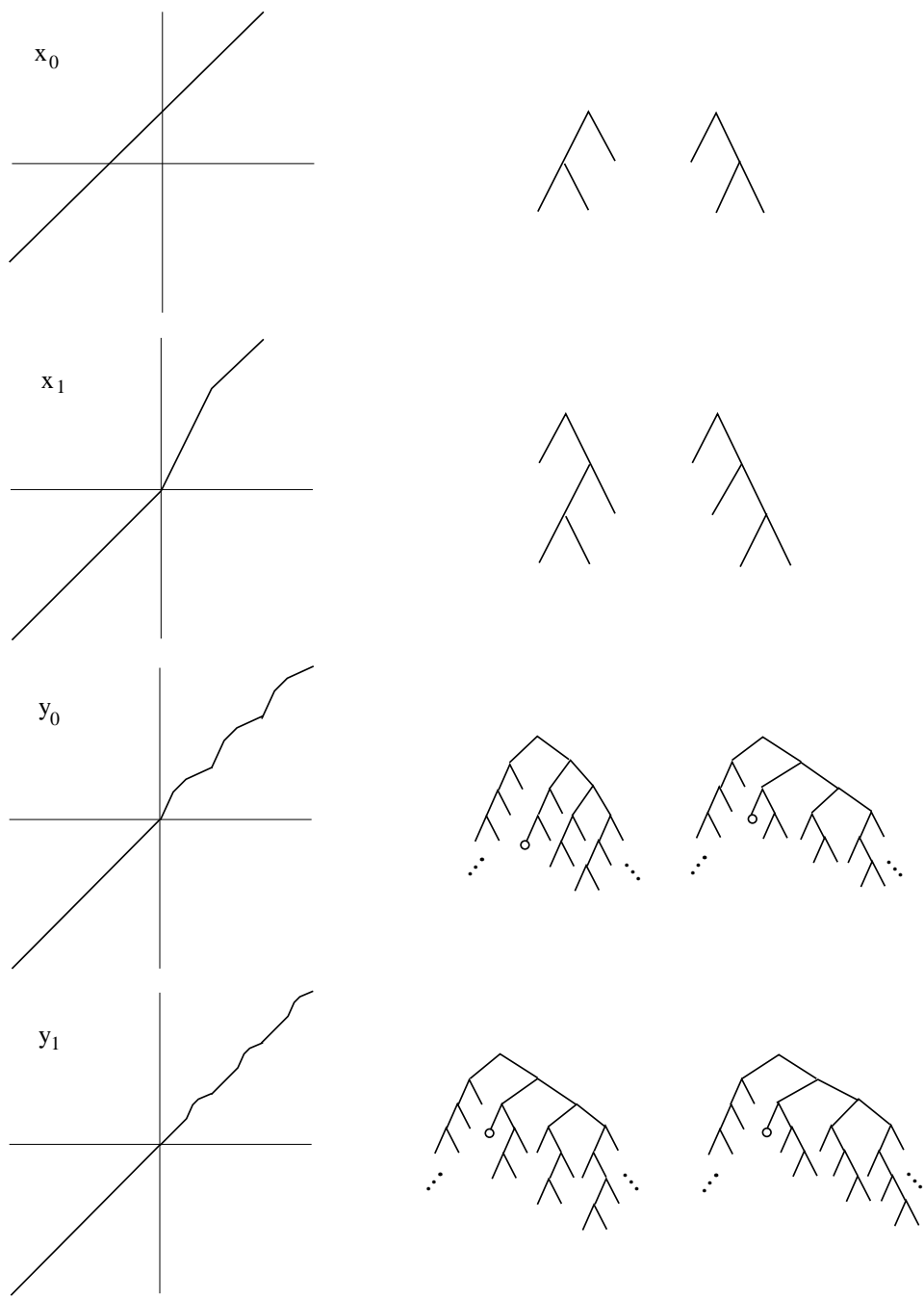}\\
  \caption{The first four generators of $\autpf$, both as maps of $\rr$ and as pairs of periodic trees. For $x_0$ and $x_1$, only the finite trees are shown as if they were in $F$, it is understood that for the periodic form, there are infinite tails at both sides of each tree (see Figure \ref{realelement}).}\label{autfgens1}
\end{figure}

\begin{figure}
  \includegraphics[width=5in]{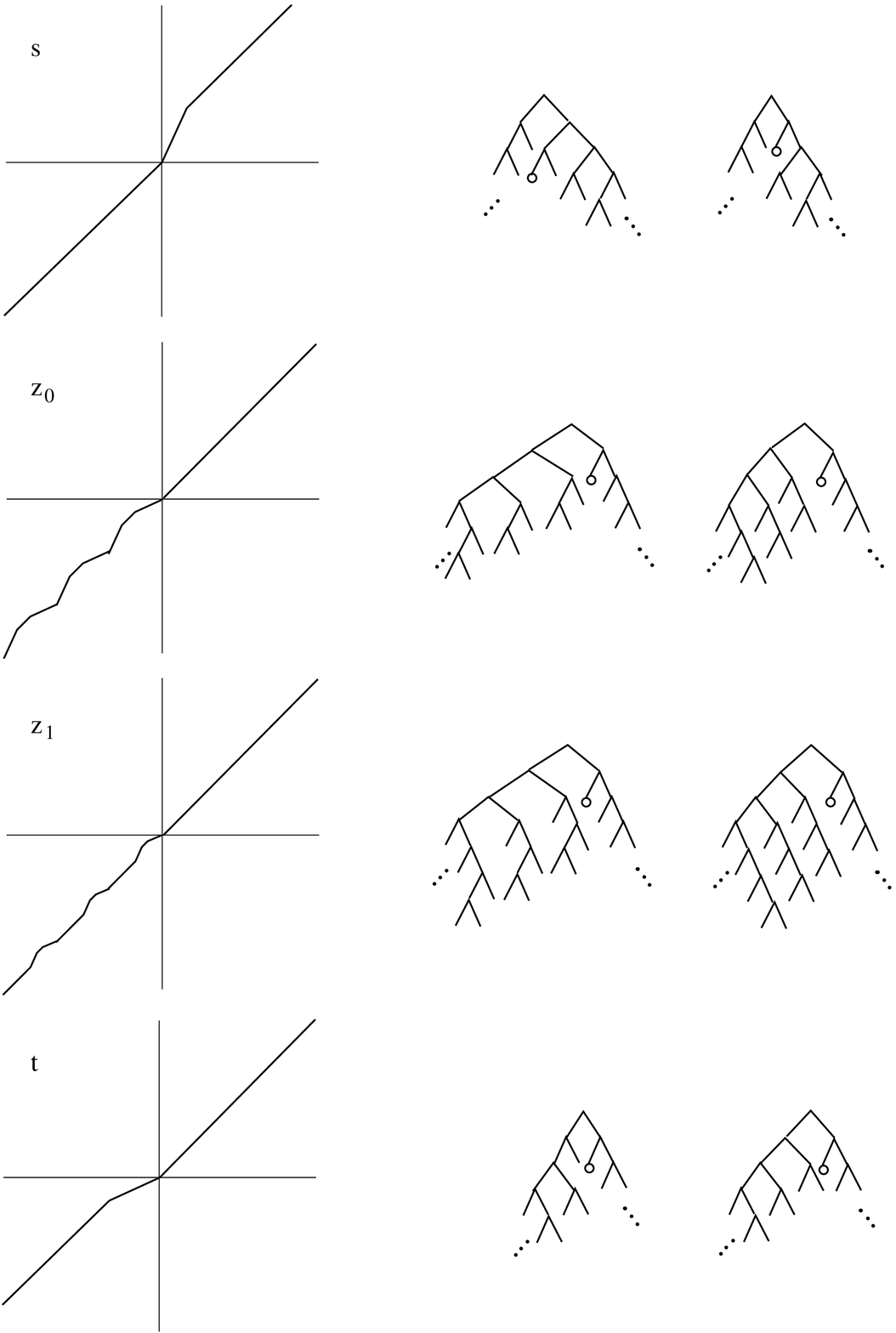}\\
  \caption{The last four generators of $\autpf$.}\label{autfgens2}
\end{figure}

Since we are using $x_0$ and $x_1$ for the $F$ embedding in $\autpf$, the two lifts of $T$ will be generated by $\{y_0,y_1,s\}$ and $\{z_0,z_1,t\}$, respectively. These lifts will be such that $y_0,y_1$ and $s$ have support in $(-\infty,0]$ whereas $z_0,z_1$ and $t$ have support in $[0,\infty)$. The lifts of $y_0,y_1,z_0$ and $z_1$ are straightforward, and amount to attaching a $[0,1]$ version of the $F$ generators into the corresponding box (from $0$ to $1$ in the case of those with support on the positive half-line, or from $-1$ to $0$ in the case of those wih support on the negative half-line), and repeating it indefinitely in the whole support. The lifts of $s$ and $t$ to $\autpf$ are trickier, and some choice is involved. Our choices are shown in Figures \ref{autfgens1} and \ref{autfgens2}, but by no means should those be considered canonical.


\begin{prop} The group $\autpf$ is presented by the generators \\$\{x_0,x_1,y_0,y_1,s,z_0,z_1,t\}$, with a set of 35 relators.
\end{prop}

The relators for $\autpf$ can be worked out following the standard procedure for short exact sequences. These relators can be classified in three types:
\begin{enumerate}
\item The two relators for the kernel $F$, namely, $[x_0x_1\ii ,x_0\ii x_1x_0]$ and $[x_0x_1\ii ,x_0^{-2}x_1x_0^2]$ (2 relators).
\item The lifts of relators for $T\times T$. These are furthermore divided into two types:
\begin{enumerate}
\item Commuting relations: he generators $y_0,y_1,s$ commute with $z_0,z_1,t$ (9 relators in total).
\item Lifts of the relators of the two copies of $T$ to $\autpf$:  Lifting the relators (f1) and (f2) is direct, they lift to the exact same relators but spelled out with $y_0,y_1$ and $z_0,z_1$. Lifting the relators (t1) to (t4) requires an extra computation and it actually depends on the choice of lifts. For instance, we will compute two examples, based on the lifts given in Figures \ref{autfgens1} and \ref{autfgens2}. The relator $s^2$ in the first copy of $T$ lifts, with our choices, to $s^2=y_1^2y_0\ii y_1\ii y_0$; and the relator (t1), in its $y_0,y_1,s$ version, lifts to
    $$
    y_1\ii sy_0^2y_1\ii y_0\ii s=y_1^4y_0^{-3}y_1y_0^2y_1y_0y_1\ii.
    $$
    There are a total of 12 relators in this class, whose tedious calculation is left to the reader.
\end{enumerate}
\item The third type of relators appearing in a short exact sequence is given by the action relators: the lifts of the generators of $T\times T$ act on the kernel $F$, since this is a normal subgroup. This gives a total of 12 relators more, each given by the conjugation action of each one of the generators $y_0,y_1,s,z_0,z_1$ or $t$ on $x_0$ or $x_1$. Again, the computation of these actions is somewhat tedious, and an example is
    $$
    sx_0s\ii=x_1^2x_0\ii x_1\ii x_0x_1\ii.
    $$
\end{enumerate}

Some of these relators will be written explicitly in the following section, when we consider subgroups of $\autpf$ generated by some of the generators.

The presentation of $\autpf$ would immediately lead to a presentation for the whole group $\autf$, again using the corresponding short exact sequence. It is worth noticing that the orientation-reversing automorphism given by $\phi(x)=-x$ would help to reduce the number of generators and relators, since the generators $z_0,z_1,t$ can be easily obtained from $y_0,y_1,t$ by conjugating by $\phi$. 
For brevity, we leave these computations to the reader as well. We have preferred to spell out the longer presentation for $\autpf$ because it will be useful to understand the subgroups studied in the next section.

\section{Some interesting subgroups}

The structure as a short exact sequence obtained for $\autpf$ shows that some interesting subgroups of it can be studied. Let
$$
\pi:\autpf\longrightarrow T\times T
$$
be the quotient map in the short exact sequence. If we consider $F$ as a subgroup of $T$, we see that we have a subgroup $F\times F$ inside $T\times T$. Clearly, these two copies of $F$ are generated by $y_0,y_1$ and $z_0,z_1$, respectively. Since we will use them often, and to clarify which one we are referring to, we will rename them $F_y$ and $F_z$. In the same way the kernel copy of $F$ will be denoted by $F_x$.

The first subgroup we will study is the group $A=\pi\ii(F_y\times F_z)$. We have now a short exact sequence, similar to the previous one:
$$
1\longrightarrow F_x\longrightarrow A\longrightarrow F_y\times F_z\longrightarrow 1.
$$
This short exact sequence has a fundamental difference with the one for $\autpf$. Each element of $F_y\times F_z$ admits a canonical lifting into $A$, since elements in $F$, seen as maps of $[0,1]$, fix 0 and 1. The elements $y_0$ and $y_1$ are lifts of the main generators of $F$ to the first component, and $z_0$ and $z_1$ to the second component. Since these elements fix the points $(k,k)$ for any $k\in \zz$, we see that these lifts define homomorphisms. Any element of $F$ can be lifted to the first or second component, by introducing its graph (in [0,1]) in the corresponding one-by-one boxes of the plane. Note that all elements obtained as lifts of elements in $F_y\times F_z$ have their graphs included in the union of all the boxes of side $[k,k+1]$, for all integers $k$.

Hence, we have proved the following result:
\begin{thm} The short exact sequence
$$
1\longrightarrow F_x\longrightarrow A\longrightarrow F_y\times F_z\longrightarrow 1
$$
splits. This gives the group $A$ a structure of semidirect product $F_x\rtimes(F_y\times F_z)$ where the action of $F_y\times F_z$ onto $F_x$ is given by the conjugation of the elements of $F_x$ by the elements $y_0,y_1$ and $z_0,z_1$.
\end{thm}

The group $A$ is the subgroup of $\autpf$ generated by $\{x_0,x_1,y_0,y_1,z_0,z_1\}$. Its 18 relators can also be classified in three types, as before:
\begin{enumerate}
\item the $F$ relators on $F_x$,
\item the $F$ relators on $F_y$ and on $F_z$, which lift unmodified to $A$ due to the splitting, and also the commuting relations: each of the $y_0,y_1$ commutes with the $z_0,z_1$,
\item the eight relators of the third type, which specify how the generators $y_0,y_1$ and $z_0,z_1$ act on the $x_0,x_1$. The first two examples are
    $$
    y_0x_0y_0\ii=x_0x_1x_0\ii  x_1x_0x_1^{-2}\qquad y_0x_1y_0\ii=x_1^3x_0\ii x_1x_0^{-3}x_1x_0x_1x_0^3x_1^{-2}
    $$
    and the rest are left to the reader. These eight relators also belong to the presentation for $\autpf$ above. As an easy example, $z_0,z_1$ commute with $x_1$ because of disjoint supports.
\end{enumerate}

Inside this subgroup $A$ we can obtain two versions of $F\rtimes F$ by just taking the first component or the second component only, giving two subgroups $B_1=F_x\rtimes F_y$ and $B_2=F_x\rtimes F_z$, both isomorphic to $F\rtimes F$ but with different actions, and whose generators are $\{x_0,x_1,y_0,y_1\}$ and $\{x_0,x_1,z_0,z_1\}$, respectively. For instance, many relators of a presentation for $B_1$ are already computed above, and we can have the actual presentation as
\begin{align*}
\<x_0,x_1,y_0,y_1\,|\,&[x_0x_1\ii ,x_0\ii x_1x_0],[x_0x_1\ii ,x_0^{-2}x_1x_0^2],\\
&[y_0y_1\ii ,y_0\ii y_1y_0],[y_0y_1\ii ,y_0^{-2}y_1y_0^2],\\
&y_0x_0y_0\ii=x_0x_1x_0\ii  x_1x_0x_1^{-2},\\
&y_0x_1y_0\ii=x_1^3x_0\ii x_1x_0^{-3}x_1x_0x_1x_0^3x_1^{-2}\\
&y_1x_0y_1\ii=x_0x_1x_0\ii x_1x_0\ii x_1x_0x_1^{-2}x_0x_1\ii\\
&y_1x_1y_1\ii=x_1^2x_0\ii x_1x_0\ii x_1x_0^{-2}x_1^2x_0^{-2}x_1x_0x_1^{-2}x_0x_1\ii x_0^3x_1^{-2}x_0x_1\ii\>
\end{align*}

However, the most interesting subgroup of $\autpf$ is obtained taking the diagonal copy of $F$ inside $F_y\times F_z$. Let $D$ be this diagonal:
$$
D=\{(g,g)\in F_y\times F_z\,|\,g\in F\}
$$
and consider $C=\pi\ii(D)$ as a subgroup of $A$, and hence, of $\autpf$. We intend to study the group $C$ extensively, since its properties are quite interesting.

The group $D$ is isomorphic to $F$, and the isomorphism is easy to understand. Take an element of $F$ as a map of $[0,1]$, and construct a map from $\rr$ to $\rr$ by introducing the one-by-one map in each box of the type $[k,k+1]\times[k,k+1]$, for each $k\in\zz$. This map is the corresponding element in $D$. Notice too that these elements act on $F_x$, seen now as maps of the real line. Hence, $D$ acts on $F_x$ by conjugation, and the result is the group $C$, a copy of $F\rtimes F$.

\begin{figure}
  \includegraphics[width=5in]{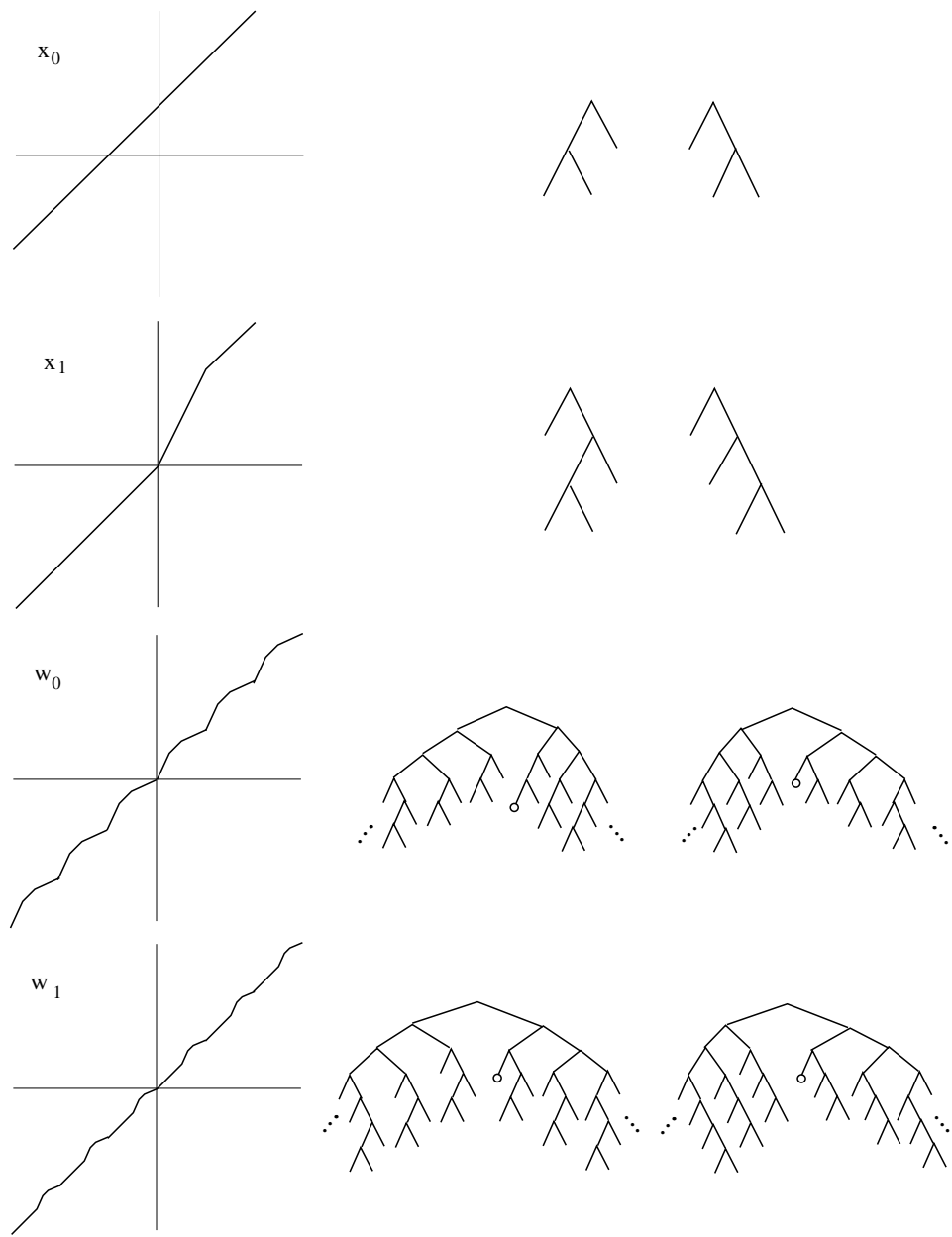}\\
  \caption{The generators for $C$.}\label{cgens}
\end{figure}

To generate $D$, since it is isomorphic to $F$, it is generated by the image of the generators for $F$. Hence, $D$ is generated by the two elements $w_0=y_0z_0$ and $w_1=y_1z_1$. For consistency, we will denote this copy of $F$ by $F_w$ as well. So the elements of $F_w$ are characterized for having the same map in each box $[k,k+1]\times[k,k+1]$, and when seen as a pair of eventually periodic trees, an element of $F_w$ has actually completely periodic trees: each one of the integral leaves (see Figure \ref{realtree}) has the same tree attached. See Figure \ref{cgens} for a picture of the generators of $C$, including the two generators for $D$.

A presentation for $C$ can be given as usual, by giving the action of $w_0,w_1$ on $x_0,x_1$. Since the maps in $F_w$ are periodic, and $x_0$ is just the map $t\mapsto t+1$ in $\rr$, then the elements of $F_w$ commute with $x_0$, something that can be seen easily with the trees as well. Hence we only need to compute the two actions of $w_0$ and $w_1$ on $x_1$. The presentation for $C$ is:
\begin{align*}
\<x_0,x_1,w_0,w_1\,|\,&[x_0x_1\ii ,x_0\ii x_1x_0],[x_0x_1\ii ,x_0^{-2}x_1x_0^2],\\
&[w_0w_1\ii ,w_0\ii w_1w_0],[w_0w_1\ii ,w_0^{-2}w_1w_0^2],\\
&[x_0,w_0],[x_0,w_1]\\
&w_0x_1w_0\ii=x_1^3x_0\ii x_1x_0^{-3}x_1x_0x_1x_0^3x_1^{-2}\\
&w_1x_1w_1\ii=x_1^2x_0\ii x_1x_0\ii x_1x_0^{-2}x_1^2x_0^{-2}x_1x_0x_1^{-2}x_0x_1\ii x_0^3x_1^{-2}x_0x_1\ii\>
\end{align*}

Note also that $w_0x_1w_0\ii=y_0x_1y_0\ii$, and $w_1x_1w_1\ii=y_1x_1y_1\ii$ (compare the presentations for $B_1$ and $C$). This is because the support of $x_1$ is $[0,\infty)$, and inside this interval, we have $w_0=y_0$ and $w_1=y_1$.

The reader is encouraged to verify these relators using periodic trees. The fact that the conjugates of elements of $F_x$ by elements of $F_w$ stay in $F_x$ remains quite apparent with the trees, since in the final element, the two tails at both ends can always be cancelled and the result trees represent an element of $F_x$.

These subgroups give new groups, all of them obtained as semidirect products of copies of $F$, and interesting actions of $F$ on itself, mixing the two interpretations of $F$ as maps of $[0,1]$ and maps of $\rr$.
In the next section we will study the metric properties of $C$, and the distortion of the two canonical embeddings of $F$ in $C$.

\section{Distortion of $F$ in $C$}\label{distortion}

The two embeddings of $F$ in $C$ behave in different ways with respect to the metric properties. Recall that in $F$, the norm of an element (the shortest word on $x_0,x_1$ representing it) is equivalent (up to a multiplicative constant) to the number of carets of the minimal reduced diagram representing it. We will use this fact to estimate the norm of an element of $F$, and compare it with its norm inside $C$. The result is the following.

\begin{thm}\label{main} The subgroup $F_x$ is quadratically distorted in $C$, whereas $F_w$ is undistorted.
\end{thm}

Recall that the distortion function of a subgroup is obtained when elements shorten their length when embedded in the large group. If $H<G$, and both are finitely generated, we have the distortion function
$$
\delta(n)=\max\{|x|_H,\text{ for all }x\in G,\ |x|_G\le n\}
$$
This theorem, stating that a copy of $F$ inside $C$ is quadratically distorted, represents an interesting change, since many of the standard subgroups of $F$ are undistorted (see Burillo \cite{burillo}, Burillo, Cleary and Stein \cite{bcs}, Cleary \cite{wreathdistort}, Cleary and Taback \cite{ctcomb}, and Guba and Sapir \cite{gubasapir}), as well as the embeddings $F\subset T\subset V$ (Burillo, Cleary, Stein and Taback \cite{thompt}). However, $F$ does admit subgroups which are at least polynomially distorted, see Guba and Sapir \cite{gubasapir}. Recently, some exponentially distorted embeddings of generalizations of Thompson's groups have appeared (see Wladis \cite{wladis}, Burillo and Cleary \cite{2v}) but they always involve more complicated groups whose elements are represented by two different types of carets. This is an example of a quadratically distorted embedding different from earlier phenomena.

The proof of Theorem \ref{main} will be spread throughout the rest of the article. We  start by proving a lower bound for the quadratic distortion of $F_x$ inside $C$ with examples.  We consider the family of words $r_n=w_0^{-n}x_1^nw_0^n$. Clearly, these words have length at most $3n$ in $C$. We will prove that each $r_n$ represents an element of $F_x$ with a number of carets of the order of $n^2$. This is enough to prove the distortion is at least quadratic.

\begin{figure}
  \includegraphics{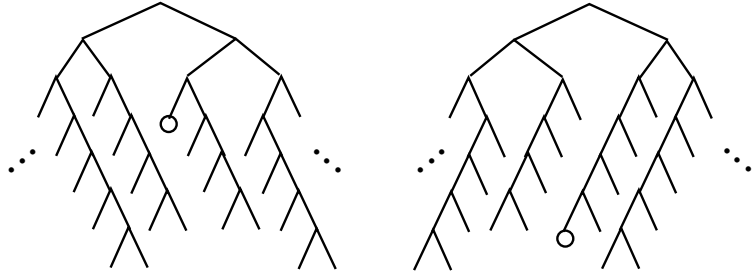}\\
  \caption{The word $w_0^4$ to start the construction of $r_4$. Notice the four-caret tree repeated throughout.}\label{r4a}
\end{figure}

\begin{figure}
  \includegraphics[width=12cm]{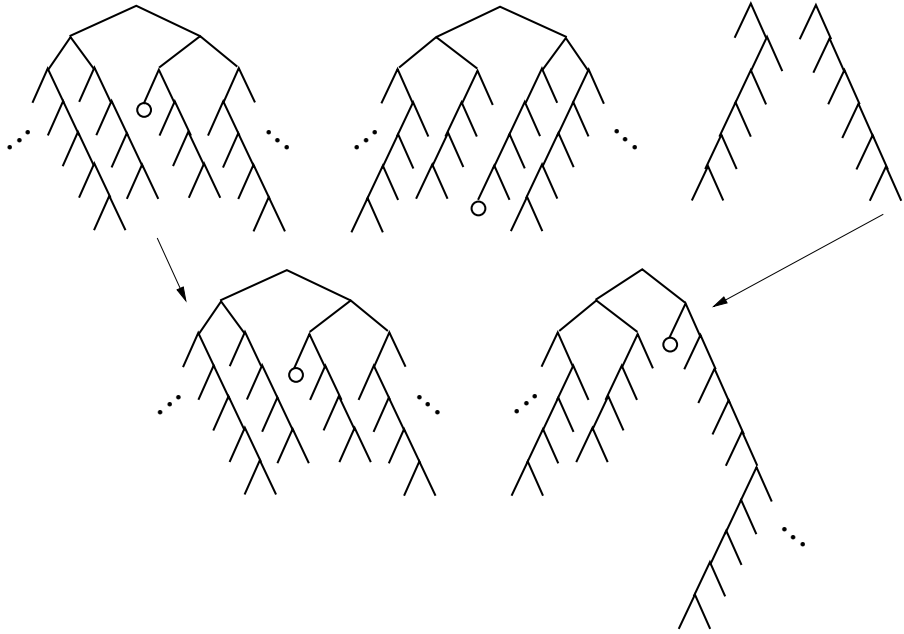}\\
  \caption{Multiplication of the word $w_0^4$ with $x_1^4$ to shift trees down the right-hand-side. Observe the empty area created, which will be filled by quadratically many carets.}\label{r4b}
\end{figure}

\begin{figure}
  \includegraphics{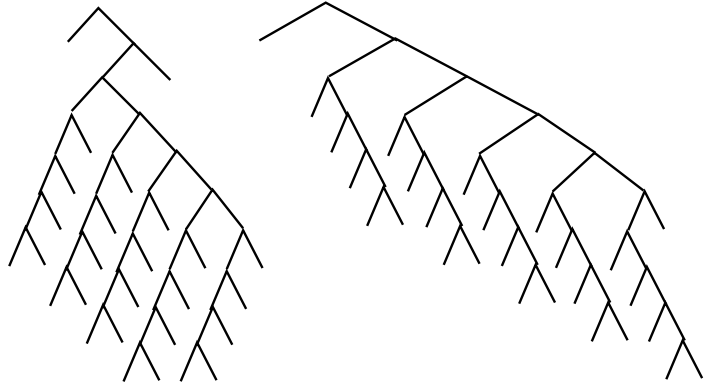}\\
  \caption{The element represented by the word $r_4$ in $F$, obtained by multiplication by $w_0^4$, which fills the five empty spaces with four-caret trees.}\label{r4c}
\end{figure}

For simplicity, we will draw pictures of a word of this type with $n=4$, but the effect will be quite clear. As shown in Figure \ref{r4a}, that the trees for $w_0^{-4}$ have a repeated tree with 4 carets in each integral leaf. Multiplying by $x_1^4$ has the effect of shifting four of those trees further down the right-hand-side of the tree, opening up a sequence of 4 empty carets on the right arm of the tree, shown in Figure \ref{r4b}. When bringing back again $w_0^4$,  those formerly empty spots are filled with trees of four carets each (see Figure \ref{r4c}), hence having at least 16 carets (actually 27). This process will give for a general $n$ a quadratically distorted word, with more than $n^2$ carets.


To prove that the distortion is no worse than quadratic, and also the linear distortion of the subgroup $\<w_0,w_1\>$, we will give a lower bound on the metric for an element of $C$ in terms of the different features of its eventually periodic tree diagram, described in the next section.

\section{Metric properties of $C$}

In this section we will relate the metric of the group $C$ with features of its elements, as it happens with $F$. Given an element of $F$, its norm with respect with the generating set $\{x_0,x_1\}$ is equivalent to the number of carets of its reduced tree diagram (see \cite{bcs} and Fordham \cite{blakegd}). Since $C$ is a semidirect product of two copies of $F$, and the elements also admit tree pair diagrams which are related to those of the copies of $F$, then it makes sense to expect similar results for the elements of $F$. We have already seen that the inclusion is distorted (at least quadratically), so a completely analogous result is not possible, but we can find a bound for the metric which guarantees that the quadratic distortion is the worst possible.

We consider the short exact sequence defining $C$
$$
1\longrightarrow F_x\longrightarrow C\stackrel\pi\longrightarrow F_w\longrightarrow 1,
$$
and let
$$
\sigma:F_w\longrightarrow C
$$
be the section of the map for $\pi$ obtained by repeating the element of $F_w$ thought of as lying in $[0,1]$ indefinitely in both directions to lie in $\rr$,  satisfying $\pi\circ\sigma=id$. Given an element $g\in C$, its periodic tree diagram is periodic with the same periodic map at both ends, and the single repeating trees give an element of $F_w$, which is exactly $\pi(g)$. The number of carets of the trees of $\pi(g)$ is the first feature we will consider, it will be called $a(g)$ in Definition \ref{metric}. This number of carets relates to the $F_w$ part of $g$. To consider the $F_x$ part of $g$, we need to see which part of $g$ is not related to the periodic parts. So we take $\pi(g)$, which is the $F_w$ part, lift it back to $C$, and use it to cancel the periodic part. This leads to the following definition.

\begin{defn} Let $g$ be an element of $C$. Then, the element of $C$ given by $d_g=(\sigma\circ\pi(g))\ii g   $ is actually an element of $F_x$, which we will call the \emph{debris} of $g$.
\end{defn}

Clearly $g$ and $\sigma\circ\pi(g)$ both map into $\pi(g)$ by $\pi$, so $d_g$ is in the kernel of $\pi$, which is $F_x$. Related to this element we will consider  two other quantities $b$ and $c$ which estimate the complexity and help us to understand the metric in $C$.

\begin{defn}\label{metric} Given an element $g\in C$, we will define the following three numbers associated to it:
\begin{enumerate}
\item The number $a(g)$ is, by definition, the number of carets of the reduced diagram of the element $\pi(g)$ in $F_w$.
\item The number $b(g)$ is the same but for $d_g$, that is, the number of carets of the reduced diagram of $d_g\in F_x$.
\item The number $c(g)$ is also computed with the reduced diagram of the element $d_g$. It is given by the sum of the number of left and right carets (not counting the root) in both trees of the reduced tree diagram for $d_g$.
\end{enumerate}
\end{defn}

\begin{figure}
  \includegraphics[width=125mm]{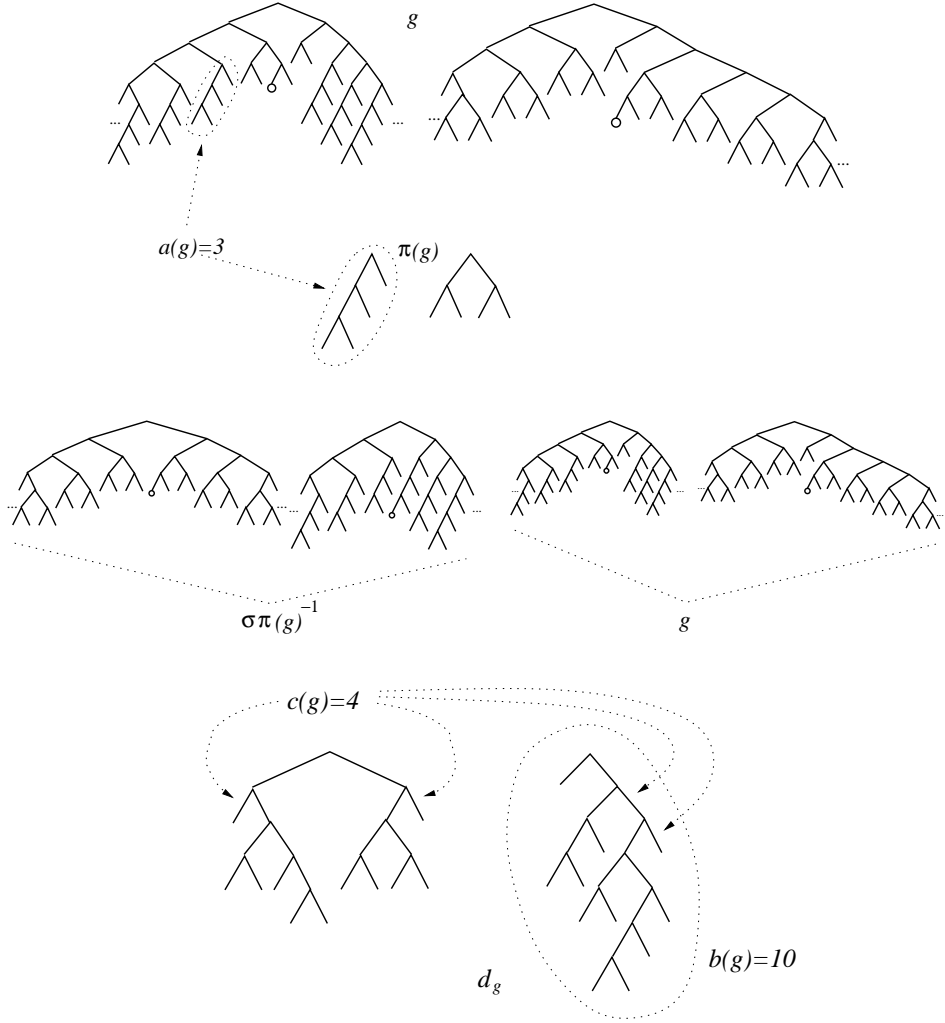}\\
  \caption{An example of an element in $C$ and its values for $a$, $b$ and $c$. From top to bottom, we see the element with $a$, the process of computing $d_g$ and the values of $b$ and $c$.}\label{abc}
\end{figure}

See Figure \ref{abc} to clarify these notions. We would like to see what is the effect on these three quantities for $g$ when it is multiplied by each one of the four generators. Studying these effects carefully we will be able to see which is the worst possible scenario for $a,b,c$ in an element of length $L$, and use it to find a lower bound for the metric.

\begin{enumerate}
\item {\bf Multiplication by $x_0$.} We see that $\pi(g)$ and $\pi(gx_0)$ are equal, corresponding to the periodic part of these elements, which is unaffected by $x_0$. From here, we see:
    $$
    d_{gx_0}=(\sigma\circ\pi(gx_0))\ii gx_0=(\sigma\circ\pi(g))\ii gx_0=d_gx_0.
    $$
    So the effect of multiplication by $x_0$ is seen only on $d_g$, and  we can consider the multiplication $d_gx_0$ as if it took place inside $F_x$.  The quantity $a$ itself will remain unchanged, whereas $b$ and $c$ may grow, but by at most two carets. So we have:
    $$
    a(gx_0)=a(g)\qquad b(gx_0)\le b(g)+2\qquad c(gx_0)\le c(g)+2.
    $$
\item {\bf Multiplication by $x_1$.} The situation is quite similar to the previous case. Multiplying by $x_1$ only affects the debris, so $a$ stays the same, $b$ and $c$ can grow by at most 3.
    $$
    a(gx_1)=a(g)\qquad b(gx_1)\le b(g)+3\qquad c(gx_1)\le c(g)+3.
    $$
\item {\bf Multiplication by $w_0$.} When we multiply $g$ by $w_0$, we notice that $a$ can grow by two carets, in case that there is no periodic part in $g$. If there is already a periodic part, it grows by at most one. We see that $c$ is unchanged from $g$ to $gw_0$. The integral leaves that will get canceled when multiplying $(\sigma\circ\pi(g))\ii$ by $g$ will get modified by $w_0$, but they will still be the same for the two elements $(\sigma\circ\pi(gw_0))\ii$ and $gw_0$, and they will get canceled the same way. So the left and right carets for $d_{gw_0}$ will be the same than those for $d_g$.

    Another way of seeing that $c$ is unchanged is by appealing to the maps in $\rr$. If $g$ is periodic outside of a box $[p,q]\times[r,s]$, for integers $p,q,r,s$, then $gw_0$ is also periodic outside this box. The element $w_0$ sends the integer interval $[r,s]$ to itself.

    Finally, if $c$ is unchanged, we see that the worst that can happen to $b$ is that it increases by two carets in every integral leaf that survives in the debris. We note that when premultiplying by  $(\sigma\circ\pi(gw_0))\ii$, in each integral leaf there may be two more carets than before, and in the debris, there could be 2 carets more surviving in each integral leaf. Hence, $b$ grows by at most $2c$. Summarizing:
    $$
    a(gw_0)\le a(g)+2\qquad b(gw_0)\le b(g)+2c(g)\qquad c(gw_0)\le c(g).
    $$
\item {\bf Multiplication by $w_1$.} The situation is again similar to the one for $w_0$. The number $a$ can grow by at most 3, $c$ stays the same, and $b$ can grow in the worst case by 3 times $c$. So we have
    $$
    a(gw_1)\le a(g)+3\qquad b(gw_1)\le b(g)+3c(g)\qquad c(gw_1)\le c(g).
    $$
\end{enumerate}

We note that multiplying by the inverses of the generators will give the same
bounds for $a$, $b$, and $c$. 

The inequalities for $b(gw_0)$ and $b(gw_1)$ where the number of additional carets can grow by a multiple of $c(g)$ is the reason of the quadratic distortion. In the example shown in Section \ref{distortion} each surviving integral leaf in $r_n$ (of which there are $n$), gets $n$ carets, giving $n^2$ carets in the element.

Armed with these inequalities, we now consider any element $g$ in $C$, and assume it has length $L=|g|_C$. So it has been constructed by multiplication of $L$ generators, and hence in each one of the $L$ multiplications, the quantities $a$, $b$ and $c$ may have suffered the worst increase detailed above. This gives the following inequalities:

\begin{enumerate}
\item $a(g)$ and $c(g)$ may have grown by at most 3 in each step, so $a(g)\le 3L$ and $c(g)\le 3L$.
\item $b(g)$ grows, in its worst case, by $3c$. Since $c$ is at most $3L$ at the end and thus also in each step, clearly we have that $b(g)\le 9L^2$, or, in a better way for our purposes, $\sqrt{b(g)}\le 3L$.
\end{enumerate}

So the conclusion is the following: there exist a constant $K>0$ such that
$$
a(g)+\sqrt{b(g)}+c(g)\le K|g|_C.
$$
This inequality is all we need to finish the proof of Theorem \ref{main}. Let $w$ be an element of $F_w$ with $N(w)$ carets, and recall that $N(w)$ is equivalent to $|w|_{F_w}$. When lifted to $C$, the element $\sigma(w)$ is all periodic, so it has $a(\sigma(w))=N(w)$, and $d_w=1$, so $b(\sigma(w))=0$ and $c(\sigma(w))=0$. Hence, $N(w)$ is bounded above by $|w|_C$ (up to a multiplicative constant), and so is $|w|_{F_w}$. Thus the periodic subgroup $F_w$ is undistorted.

For the embedding of the nonperiodic subgroup $F_x$, we take $x$ a general element of $F_x$ and consider it also as an element of $C$. When seen in $F_x$, it has some number of carets $N(x)$, equivalent to $|x|_{F_x}$. But naturally $d_x=x$ if $x\in F_x$, so we have $b(x)=N(x)$. So our inequalities above say that there exists a constant $K$ such that
$$
|x|_{F_x}\le K|x|_C^2
$$
which proves that the distortion of $F_x$ in $C$ is at most quadratic. From the examples explicitly constructed in Section \ref{distortion}, we see that the distortion is exactly quadratic.

\bibliographystyle{plain}

\def\cprime{$'$} \def\cprime{$'$}

\end{document}